
Made for Amstex 2.1; inputs amsppt.sty

\documentstyle{amsppt}
\magnification=\magstep1
\nologo

\def\proof{\demo{Proof}}
\def\endproof{\qed\enddemo}
\def\nextpart#1{\medskip{\bf (#1).}}

\def\set{\!:=}
\def\nmb{\nomathbreak}

\def\tl#1{\widetilde{#1}}
\def\oY{{\overline{Y\mkern-.5mu}}}

\def\ord{{\text{\rm ord}}}
\def\Spec{\text{\rm Spec}}
\def\Proj{\text{\rm Proj}}
\def\Hom{{\text{\rm Hom}}}
\def\sHom{{\Cal H\eurm o\eurm m}\mkern2mu}
\def\I{{\Cal I}}
\def\J{{\Cal J}}
\def\O{{\Cal O}}
\def\C{{\frak C}}

\loadeurm
\loadbold

\topmatter

\title {Adjoints  of ideals in regular local rings}
\endtitle

\author Joseph Lipman \endauthor

\affil Purdue University \endaffil

\thanks Partially supported by the National Security
Agency.\vadjust{\kern1.3pt}  \endthanks

\address{Department of Mathematics, Purdue University,
 W. Lafayette, IN 47907, USA} \endaddress

\email{lipman\@math.purdue.edu}\endemail

\subjclass  13H05, 13C99, 14F17 \endsubjclass

\endtopmatter

\document

\subheading{Introduction}
Several existing results labeled ``Brian\c con-Skoda theorem'' concern
an ideal~$I$ in a regular local ring~$R$. Of these, the weakest states
that if $I$ is generated by $\ell$ elements, then $\overline {I^{n+\ell
-1}}\subset I^n\ (n>0)$, where ``$\mkern2.3mu\overline{\phantom{I}}\,$" denotes
``\kern.4pt integral closure."
In this paper, we associate to $I$ an integrally closed ideal $\tl
I\supset\bar I$, the  {\it adjoint\/} of ~$I$, and indicate how it can
be used in place of~$\bar I$ to improve such results. At first, in
Theorem~(1.4.1), this just involves a recycling of methods
from~\cite{LS}. (Even that is not without benefit, see Corollary~
(1.4)).  But there's more. It's not hard to show that there is an $n_0$
such that $\tl{I^{n+1}}=I\mkern1mu\tl{I^n}$ for all $n\ge n_0$.  The basic
conjecture~(1.6)---which, as we'll see, quickly implies a number of
recently proved Brian\c con-Skoda-type theorems---says that $n_0$ can
be taken to be $\ell(I\mkern1mu)-1$, where $\ell(I\mkern1mu)$ is the analytic
spread
of~$I$.  This conjecture does hold when $R$ is essentially of finite
type over a field of characteristic zero, or when $\dim R=2$.

 Section 2 deals with a conjecture, related to  Grauert-Riemenschneider
vanishing, about certain cohomology groups being zero.
Suppose there exists a proper birational map $f\:Y\to\nmb \Spec(R)$ such
that $I\O_Y$ is invertible and $Y$ is {\it nonsingular,\/} i.e.,
locally regular.  (The existence of such a $Y$ in all characteristics
is not yet certain, but it is needed in the vanishing conjecture.) Let
$\omega_Y$ be a dualizing sheaf for ~$f$, chosen to be canonical in the
sense that its restriction to the open set~$U$ where $f$ is an
isomorphism is $\O_U$. While the definition (1.1) of $\tl I$ uses
neither $Y$ nor any duality theory, Proposition~(1.3.1) states that
$$
\tl I\set H^0(Y, I\omega_Y\mkern-2mu)\mkern1mu;
$$
and the vanishing conjecture states that
$$
H^i(Y, I\omega_Y\mkern-2mu)=0\qquad(i>0).
$$
The point is that this vanishing conjecture implies conjecture (1.6).
In fact it is thanks to Cutkosky's transcendental proof of the
vanishing conjecture \cite{C} that we know (1.6) holds in
characteristic zero.

Section 3 elaborates on the two-dimensional case (where the
vanishing conjecture is known to hold, see Remark (2.2.1)(b)).
A geometrically motivated treatment
of the adjoint of a simple complete ideal $I$  is given in~\cite{L4};
close connections with the multiplicity sequence and the
conductor ideal of the local ring~$\frak o$
of the ``generic curve through ~$I\mkern1.5mu$" are brought out. Roughly
speaking,
{}~$\frak o$ is the local ring at the generic point of the exceptional divisor
(a $\Bbb P^1$) on the blowup~$Y_0$ of any 2-generated reduction ~$I_0$ of~$I$.
Propositions~(3.1.1) and (3.1.2) below explore such connections for an
arbitrary integrally closed~$I$.  If $Y$ is the normalization of~$Y_0$,
then $I\omega_Y$ is just the conductor~$\C\set\O_{Y_0}\mkern-1mu\:\!\O_Y$
(so $\C$ is independent of the choice of~$I_0$), and it is generated by its
global sections~$\tl I$. We also find in~Proposition~ (3.3) that
$\tl I=I_0\:\!I$.

Furthermore, in~ \cite{HS} Huneke and Swanson have shown that $\tl I$
is just the second Fitting ideal $\frak F_2(I\mkern1mu)$; and more
generally (since $I=\frak F_1(I\mkern1mu)$), that for all $n>0$, $\;\tl
{\frak F_n(I\mkern1mu)}=\frak F_{n+1}(I\mkern1mu)$.

Let me mention in closing that though the material in~\cite{L4} dates
back to 1966, the results in this paper all came out of an effort to
analyze the Brian\c con-Skoda theorem~(3.3) in~\cite{AH1}.

\bigskip

\subheading{1. Adjoints and Brian\c con-Skoda theorems}
Let $R$ be a regular noetherian domain with fraction field~$K\mkern-1.5mu$,
let $v$ be a valuation of~$K$ whose valuation ring ~$R_v$
(with  maximal ideal $m_v$) contains~$R$,
and let~$h$ be the height of the prime ideal $p\set m_v\cap R$.
We say that
$v$ is a {\it prime \-divisor\/} of~$R$ if $R_v/m_v$ has transcendence
degree~$h-\nmb1$ over its subfield~$R_p/pR_p\mkern1mu$. It is equivalent that
$R_v$~be essentially of finite type over~$R$, or that $v$
be a Rees valuation of some $R$-ideal ~$I\mkern-1mu$, i.e., that $R_v$ be
$R$-isomorphic to the local ring of a point on the normalized blowup
$\oY_{\!I}\set\Proj(\oplus_{n\ge 0}\,\overline{I^n}\mkern1mu)$,
where $\overline{I^n}$ is the integral closure
of~$I^n$. Such a~$v$ is a discrete rank-one valuation.
(See \cite{A, p.\,300, Thm.\,1\,(4) and p.\,336, Prop.\,3}. Note also that
$R$, being universally catenary, satisfies the ``dimension formula"
\cite{EGA\,III, (5.6.4) and (5.6.1)\,(c)};
and that $\oY_{\!I}$ is of finite type over~$R$ \cite{R, p.\,27, Thm.\,1.5}.)

\definition{Definition (1.1)} The {\it adjoint\/} of an $R$-ideal~$I$ is the
ideal
$$
\tl I\set \bigcap_v\,\{\,r\in K\mid v(r)\ge v(I\mkern1mu) - v(J_{R_v/R})\,\}
$$
where the intersection is taken over all prime divisors~$v$ of~$R$,
\footnote{We consider two valuations with the same valuation ring to be
identical.}
and for any essentially finite-type $R$-subalgebra~$S$ of~$K\mkern-1.3mu$,
the {\it Jacobian ideal\/}~$J_{S/R}$ is the 0-th Fitting ideal
of the $S$-module of K\"ahler differentials~$\Omega^1_{S/R}\mkern1.5mu$.
\enddefinition

\remark\nofrills{Remarks {\bf(1.2).}\usualspace} (a) $\tl I\subset R$ because
$R$ is the intersection of its localizations at height one primes,
and each such localization is the
valuation ring of a~$v$ for which $v(J_{R_v/R})= 0$. Hence
$$
\tl I= \bigcap_v\,\{\,r\in R\mid v(r)\ge v(I\mkern1mu) - v(J_{R_v/R})\,\}
$$
where the intersection is taken over all prime divisors~$v$\vadjust{\kern1pt}
such that $v(I\mkern1mu)>0$.

(b) Being an intersection of valuation ideals,
$\tl I$ is {\it integrally closed;}\vadjust{\kern1pt}
and if $\bar I$ is the integral closure of~$I$ then
$$
I\subset \bar I \subset \tl I = \,@,@,\tl{\!@!@!\bar I}.
$$

(c) For any $x\in R$, we have $\tl{xI}=x\tl I$. In particular, $\tl{xR}=xR$.
\vadjust{\kern1pt}

(d) For any two $R$-ideals $I$, $J$, we have $\tl{JI}\:\! I=\tl J$.
In particular,
$$
\tl{I^{n+1}}:I=\tl{I^n} \qquad(n\ge 0).
$$
\endremark

\nextpart{1.3}
For any finite-type birational map $f\:Y\to\Spec(R)$,
we may---and will---identify $\O_Y$ with a subsheaf of the constant sheaf~$K$
on~$Y\mkern-1mu$, so that the stalks~$\O_{Y\!,\mkern1.3mu y}\ (y\in Y)$
are all $R$-subalgebras of~$K$.
If $g\:Z\to\Spec(R)$ is another such map which factors via ~$f$, then $g$ is
uniquely determined by $Z$ and $Y\mkern-2mu$, and we say that $Z$ {\it
dominates\/} $Y\mkern-2mu$.
The {\it relative Jacobian\/} $\J_{\mkern-2.3mu f}$
(or, less precisely, $\J_Y$) is
the coherent $\O_Y$-module whose sections over any affine open
$\Spec(A)\subset Y$ are given by
$$
H^0\bigl(\Spec(A),\J_Y\bigr)=J_{A/R}\,.
$$
We set
$$
\omega_Y\set \O_Y\:\mkern-1mu\J_Y\cong \sHom_Y(\J_Y,\mkern1mu\O_Y).
$$
If $Y$ is normal, $\omega_Y$ is a canonical
{\it relative dualizing sheaf\/} for ~$f$ \cite{LS, p.\,206, (2.3)}.

For any {\it proper\/} birational $f\:Y\to\Spec(R)$ with $Y$
{\it normal\/} and $I\O_Y$ {\it invertible,} we set
$$
\tl I_Y \set H^0(Y, I\omega_Y),
$$
the ideal obtained by restricting the intersection in\vadjust{\kern .7pt}
Definition~(1.1)
to those ~$v$ such that $R_v$ is $\O_{Y\!,\mkern 1.3mu y}$ for some $y\in
Y\mkern-1mu$. So $\tl I\subset \tl I_Y\mkern1mu$,
and $\tl I_Y$ is a ``decreasing" function of~$Y$ in the sense that for any
proper birational $g\:Z\to Y$ with $Z$ normal, we have
$\tl I_Z\subset \tl I_Y$. For any prime divisor~$w$, $R_w$ is the local ring
of a point on some such ~$Z$; so the intersection of all
$\tl I_Z$ is just $\tl I$.

\proclaim{Proposition (1.3.1)}
For any\/ $@,@,Y@!$ as above and having pseudo\kern.5pt-rational singularities
(for example, $Y$ regular), $\tl I_Y=\tl I$.
If such a\/~$Y$ exists then for any multiplicative system~$M$
in\/~$R,$\ $\tl{IR_M}=\tl IR_M$.
\endproclaim

\proof
The pseudo\kern.3pt-rationality assumption forces $g_*(I\omega_Z)=I\omega_Y$
for all
$g\:Z\to Y$ as above (by \cite{LT, p.\,107, Corollary}, and since
$I\O_Y$ is invertible), whence
$$
\tl I_Z= H^0(Z,I\omega_Z)=H^0(Y,\, g_*I\omega_Z)=\tl I_Y,
$$
and $\tl I_Y=\cap_Z\,\tl I_Z=\tl I$.
The rest follows from compatibility of $H^0(Y,I\omega_Y)$ with
localization on ~$R$.
\endproof

\remark\nofrills{Remarks {\bf(1.3.2).}\usualspace}
(a) That a regular~$Y$ with
$I\O_Y$ invertible {\it always\/} exists has been announced by Spivakovsky,
but details have not appeared at the time of this writing. For the
equicharacteristic zero case, see \cite{H}.

(b) In dimension~2, every normal
{}~$Y$ birationally dominating~$\Spec (R)$ has pseudo\kern.3pt-rational
singularities,
\cite{L1, p.\,212, \S9}, \cite{LT, p.\,103, Example\;(a)}. So
in~Prop.\,(1.3.1),
we could take $Y$ to be the normalized blowup of~$I$.
\endremark

(c) Another case in which we could take $Y$ in~(1.3.1)
to be the blowup of~$I$ is when $I=(x_1,\dots,x_r,y)$ where
$(x_1,\dots,x_r)$ is a regular sequence such that $R/(x_1,\dots,x_r)R$ is
still regular \cite{LS, p.\,219, Proposition, (ii)}. Here
$\tl {I^n}_Y\ (n\ge\nmb 0)$ is easily calculated: indeed,
if $L$ is any $R$-ideal generated
by a regular sequence of length ~$\ell$, and such that all the powers of $L$
are integrally closed, then the blowup~$X$ of~$L$ is normal
and $\J_X=(L^{\ell-1}\O_X)$,\,
\footnote
{because for a regular sequence $(a_1,\dots,a_\ell)$, we have, e.g.,
$$
R_\ell\set R[\frac{a_1}{a_\ell},\dots,\frac{a_{\ell-1}}{a_\ell}]\cong
R[t_1,\dots,t_{\ell-1}]/(a_\ell t_1-a_1,\dots,a_\ell t_{\ell-1}-a_{\ell-1}),
$$
and so $\Omega^1_{R_\ell/R}$ is generated by the differentials $d(a_i/a_\ell)$
subject to $a_\ell d(a_i/a_\ell)=0\ \,(1\le\nmb i\le\nmb \ell-1)$.%
}
\,so that
$$
\alignat2
\tl {L^n}_X&=R \qquad&& (n<\ell)\\
&=L^{n-\ell+1}\qquad&&(n\ge \ell).
\endalignat
$$
For $L=I$, we have $\ell=r+1$ or $\mkern1mu r$.

\nextpart{1.4}
The following  is clearly related to the Brian\c con-Skoda theorem in
\cite{LS, p.\,204, Thm.~$1''$}.
(Recall the above-given inclusion $\bar I\subset\tl I\mkern1mu$,
where $\bar I$ and $\tl I$ are the integral closure
and adjoint, respectively, of the $R$-ideal ~$I$.)

We say that $I$ is {\it $\ell$-generated\/} ($\ell\ge 0)$ if $I$ is generated
by $\ell$ elements.

\proclaim{Theorem (1.4.1)}
For any\/ $\ell$-generated ideal\vadjust{\kern1pt} $I$
in a regular noetherian domain~$R$:

{\rm (i)} $\tl{I^{n+\ell -1}}\subset I^n$ for all\/ $n\gg0$.

{\rm (ii)} If the graded ring
$\text{gr}_IR\set\oplus_{n\ge0}\mkern2mu I^n\mkern-1mu/I^{n+1}$
contains a homogeneous regular element of positive degree, then\/
{\rm(i)} holds for all\/ $n\ge1$.

\nopagebreak
{\rm (iii)} $\tl{I^{n+\ell}}\subset I^n$ for all\/ $n\ge0$.

\endproclaim

\proof
If $Y_0\set\Proj(\oplus_{n\ge0}\mkern2mu I^n)$ is the blowup of $I$,
and $Y$ is its normalization,\vadjust{\kern.5pt} then as in \cite{LS, p.\,200,
Thm.\,2 and proof of Corollary}, we have
$I^{\ell-1}\omega_Y\subset\O_{Y_0}$
(all inside the constant sheaf~$K$ on~$Y_0$). Hence
$$
\tl {I^{n+\ell-1}}\subset
   H^0(Y, \mkern1mu I^{n+\ell-1}\omega_Y)
    \subset H^0(Y_0, I^n\O_{Y_0})
      =\bigcup_{j\ge 0}\,I^{n+j}\:\mkern-1.5mu I^j\mkern1mu.
$$

For $n\gg0$, $H^0(Y_0, I^n\O_{Y_0})=I^n$ (by e.g., \cite{EGA\,III, (2.3.1)}),
\footnote{For a more elementary proof, apply \cite{ZS, pp.\medspace
154--155, Lemmas 4 and 5} to the ideal $\frak B\set (0)$ in~$\text{gr}_IR$
to find an integer $q$ such that for any $n$ and any $x\in H^0(Y_0,
I^n\O_{Y_0})\setminus I^n$, $x\notin I^q$ (because the leading form of~$x$
annihilates all homogeneous elements of large degree \dots\!) Such an $x$
must lie in $\mkern1mu\overline {\mkern-1mu I^n}$. But there exists $p$
such that
$\mkern1mu\overline{\mkern-1mu I^{p+1}}=I\mkern2.5mu\overline{\mkern-1mu I^p}$,
whence $\mkern1mu\overline{\mkern-1mu I^{p+q}}\subset I^q$,
and therefore $n<p+q$.\looseness=-1%
}
proving~(i).

If $\text{gr}_IR$ has a homogeneous regular element of
positive degree, then $I^{n+j}\:@!@!@!I^j=\nmb I^n\!$, proving~(ii).
\smallskip
\eightpoint
In (i), the restriction of $n$ to sufficiently large values is annoying,
and may well be unnecessary (see Conjecture~(1.6) below). If so,
then (iii)---and the following ungainly argument---would be superfluous.

The polynomial ring~$R[t]$ is still regular.
An immediate consequence of the following Lemma is that
 {\it for any $R$-ideal\/~$L$,
 $\tl LR[t]\subset \tl{LR[t]}$.} (The adjoints are taken in~$R$ and~$R[t]$
respectively.) With $I'\set(I,t)R[t]$, we have that $\text{gr}_{I'}R[t]\cong
(\text{gr}_{I'}R)[t]$ has a
regular element (namely ~$t$) of degree~1, and we can apply (ii) to get
$\tl{I'{}^{n+\ell}}\subset\nmb I'{}^n$ for all\/ $n\ge0$; and since
$\tl{I^{n+\ell}}R[t]\subset(I^{n+\ell}R[t])\tl{\kern10pt}\subset \tl
{I'{}^{n+\ell}}$ and $I'{}^n\cap R=I^n$, therefore (iii) results.

\proclaim{Lemma (1.4.2)}
Let\/ $w$ be a prime divisor of the polynomial ring\/~$R[t]$ and let\/~$v$
be the restriction of~$w$ to~$K$, the fraction field of~$R$. Then\/
$v$ is a prime divisor of\/~$R$, and for any\/ $R$-ideal\/~$L$,
$$
v(L) - v(J_{R_v/R})\ge w(L) - w(J_{R_w/R[t]}).
$$
\endproclaim
\demo{Proof}
Let $(R_w,m_w)$ and $(R_v,m_v)$ be the (discrete) valuation rings of~$w$
and~$v$ respectively. Set $q\set m_w\cap R_v[t]$. There are two
cases to consider.

(1) $q=m_vR_v[t]$. Then the localization~$R_v[t]_q$ is a discrete valuation
ring contained in, and hence equal to, $R_w$. Thus $R_w/m_w=(R_v/m_v)(t)$ has
transcendence degree (t.d.)~1 over $R_v/m_v$.

(2) $q\supsetneq m_vR_v[t]$, whence $q$ is maximal, of height~2, and $R_v[t]/q$
is algebraic over~$R_v/m_v$.
Since $R_w$ is essentially of finite type over~$R[t]$,
hence over~$R_v[t]$, therefore $R_w$ is a prime divisor of ~$R_v[t]$; and
so $R_w/m_w$ has t.d.~1 over~$R_v[t]/q$. Thus, again,
$R_w/m_w$ has t.d.~1 over $R_v/m_v$.
\smallskip

Now set $p'\set m_w\cap R[t]$ and $p\set p'\cap R=m_v\cap R$, so that
$R_w/m_w$ has t.d.~$\text{height}(p')-1$ over ~$R[t]/p'$,
and, by the preceding remarks, the t.d. of ~$R_v/m_v$
over ~$R/p$ is $\text{height}(p')-2+\text{the}$
t.d\. of~$R[t]/p'$ over~$R/p$.
It follows then from \cite{ZS, p.\,323, Prop.\,1A} that $R_v/m_v$
has t.d.~$\text{height}(p)-1$ over ~$R/p$, and so $v$
is indeed a prime divisor of~$R$.

The last assertion follows from the relation
$$
J_{R_w/R[t]}=J_{R_w/R_v[t]}J_{R_v[t]/R[t]}=J_{R_w/R_v[t]}J_{R_v/R}.
$$
(See \cite{LS, p.\,201, (1.1)} for the first equality.)
\endproof
\tenpoint\enddemo

Suppose now that $R$ is {\it local,\/} with maximal ideal~$m$.
For an $R$-ideal~$I$, the {\it analytic spread\/} $\ell(I\mkern1mu)$
is the dimension of the ring $\oplus_{n\ge0}\mkern2mu I^n\mkern-1mu/mI^n$.
When $R/m$ is infinite, $I$ has an $\ell(I\mkern1mu)$-generated {\it
reduction\/}
{}~$I_0\subset I$, i.e., $I_0I^n=I^{n+1}$ for some~$n\ge 0$.

\proclaim{Corollary (1.4.3)}
For $R$ local, assertions\/ {\rm(i)} and\/
{\rm(iii)} in~Theorem\/ $(1.4.1)$ hold with\/
$\ell$~the analytic spread of~$I$. And if $I$ has an\/ $\ell$-generated
reduction~$I_0$ such that $\text{gr}_{I_0}R$
contains a homogeneous regular element of positive degree, then\/
{\rm(i)} holds for all\/ $n\ge0$.\looseness =-1

\endproclaim

\demo{Proof} By arguing as in the proof of~(1.4.1)(iii), with $R[t]$
replaced by its localization $S\set R[t]_{mR[t]}$, and $I'\set IS$, we
reduce to the case where $R/m$ is infinite. Then\vadjust{\kern.7pt} we can
apply~(1.4.1) to
an $\ell$-generated reduction~$I_0$, noting that for any valuation~$v$
such that $R_v$ contains~$R$ we have $v(I_0)=\nmb v(I\mkern1mu)$, whence
$\!\tl{\,I_0^p}=\!\tl{\,I^p}$ for all $p\ge 0$.
\endproof
\goodbreak
The following statement was conjectured by Huneke.

\proclaim{Corollary (1.4.4)}
If $(R,m)$ is a $d$-dimensional regular local ring and
$I$ is an $m$-primary ideal, then for all~$n\ge1,$
$$
\mkern1mu\overline {\mkern-1mu I^{n+d-1}}\:\!m^{d-1}\subset I^n.
$$
\endproclaim

\proof
Replacing $(R,I\mkern1mu)$ by $(S\set R[t]_{mR[t]}, \mkern1mu IS)$ if
necessary, we may
assume that $R/m$ is infinite. Then $I$ has a $d$-generated
reduction~$I_0$ such that $\text{gr}_{I_0}R$ is a polynomial ring in
$d$ variables over~$R/I_0$; so Corollary~(1.4.3) gives
$\!\tl{\,I^{n+d-1}}\subset I^n$. Thus it suffices to show that
$\mkern1mu\overline{\mkern-1mu I^{n+d-1}}\:m^{d-1}\subset\! \tl{\,I^{n+d-1}}$,
for which it's clearly enough (see (1.2)\,(a)) that
{\it for any prime divisor\/~$v$ of\/ $R$ such that\/ $m_v\cap R=m,$}
$$
v(J_{R_v/R})\ge v(m^{d-1}).
$$
But $R_v$ contains $R'\set R[x_2/x_1,\dots,x_d/x_1]$ for some generating set
$(x_1,x_2,\dots,x_d)$ of~$m$, and then
$$
J_{R_v/R}=J_{R_v/R'}J_{R'/R}=m^{d-1}J_{R_v/R'}
$$
(see \cite{LS, p.\,201, (1.1) and top of p.\,202}), which gives
the desired result.
\endproof

\proclaim{Lemma (1.5)}
Let\/ $R$ be a regular noetherian domain, let\/ $I$ be an\/ $R$-ideal, and
set\/ $G\set\nmb\oplus_{n\ge0}\,I^n,$\ $\tl{G}\set\oplus_{n\ge0}\,\tl{I^n}$.
Then\/ $\tl{G}$ is a finitely generated graded\/ $G$-module, and hence there
is an $n_0$ such that
$$
\tl{I^{n+1}}=I@,@,\tl{I^n}\quad \text{for all }\medspace n\ge n_0.
$$
\endproclaim

\proof
$\tl G$ is a graded $G$-module because, clearly,
$I^p\tl{I^q\,}\!\subset \!\tl{\,I^{p+q}}\ (p,q\ge 0)$.
Now just observe, with $Y$ the normalized blowup of~$I$, that by~(1.3),
$\tl{G}$ is a submodule of
$\oplus_{n\ge0}\,H^0(Y,I^n\omega_Y)$,
which is finitely generated over~$G$ \cite{EGA\,III, (3.3.2)}.
\endproof

As we'll see in (2.3) below, the following refinement of Lemma~(1.5)
holds true when $R$ is essentially of finite type over a
characteristic-zero field, or when $\dim R=\nmb2$.
(The 2-dimensional case also results from
Prop.\,(3.1.2); see also Prop.\,(4.2)
of \cite{HS}. For another example, see Remark (1.3.2)(c).))

\proclaim{Conjecture (1.6)}
Let $R$ be a regular local ring, and let $I$ be an $R$-ideal of
analytic spread~$\ell$. Then
$$
\tl{I^{n+1}}=I@,@,\tl{I^n}\quad\text{for all }\medspace n\ge\ell-1.
$$
\endproclaim
\smallskip
We illustrate the usefulness of this conjecture (when it holds)
by indicating how it implies some Brian\c con-Skoda-type
theorems recently proved for equicharacteristic regular local rings
by Aberbach and Huneke. These
theorems are all of the form\looseness=-1
$$
\overline{I^{n+\ell-1}}\subset I^nA\qquad(n>0),
$$
where the ``coefficient ideal"~$A$ depends only on~$I$. Under the assumption
that (1.6) holds, we need only show that $\tl{I^{\ell-1}}\subset\nmb A$ in
order
to get the stronger assertion
$$
\tl{I^{n+\ell-1}}=I^n\tl{I^{\ell-1}}\subset I^nA.
$$

\nextpart{1.6.1} In \cite{AH2}
$A$ is taken to be the sum of all ideals $A'$ such that $IA'=\bar IA'$.
By (1.6) and (1.2)(b), $I\tl{I^{\ell-1}}=\tl{I^{\ell}}=\bar I\tl{I^{\ell-1}}$,
so that $\tl{I^{\ell-1}}\subset\nmb A$.
\footnote{When $\dim R=2$ and $I$ is a 2-generated ideal primary for the
maximal ideal, then $\tl I=A$, see Prop.\,(3.3) below.}

\nextpart{1.6.2} In \cite{AH1, p.\,350, Thm.\,3.3}, $A$ is taken to be
the intersection of the primary components of~$I^{\ell-h}$ belonging to the
minimal primes $p_1,\dots,p_e$ of~$I$, where
$h\set\nmb\max_i h_i\set\nmb \max_i \text{height}(p_i)$.\vadjust{\kern1pt}
(To check that $\ell\ge h$, just localize at each ~$p_i$.) To~show
that $\tl{I^{\ell-1}}\subset\nmb A$, localize at $p=p_i\ (1\le i\le e)$,
and note that
$$
\tl {I^{\ell-1}}R_p\subset\tl{I^{\ell-1}R_p}\subset I^{\ell-h_i}
\subset I^{\ell-h},
$$
where the first inclusion is elementary, and the second is given
by~(1.4.1)(ii).

Moreover, if (1.6) holds, then, with $I_p\set IR_p$, we have
$$
\tl{I_p^{\ell-1}}=I_p^{\ell-h}I_p^{h-h_i}\tl{I_p^{h_i-1}}
 =I_p^{\ell-h}\tl{I_p^{h-1}}\,;
$$
and hence if $\tl {I^{\ell-1}}R_p=\tl{I_p^{\ell-1}}$ for all~$p_i$ (see
Prop.\,(1.3.1)),
then $\tl{I^{\ell-1}}$ is contained in the intersection of the primary
components of~$I^{\ell-h}\tl{I^{h-1}}$ belonging to the $p_i$.

\nextpart{1.6.3} In \cite{AHT, Thm.~7.6}, the above-mentioned Theorem~3.3
of \cite{AH1} is strengthened. Here the inductive description of~$A$ is
somewhat complicated. So suffice it to say that the inclusion
$\tl{I^{\ell-1}}\subset\nmb A$ can be established by
alternately localizing at suitable associated primes of height~$i$ and applying
{}~(1.6), as $i$ goes, one step at a time,
from $\ell-1$ down to the height of~$I$.

\medskip

\subheading{2. A vanishing conjecture}
Again, let~$I$ be an ideal in a regular local ring~$(R,m)$.
Throughout this section we make the following assumption---which is
satisfied at least over varieties in characteristic zero \cite{H,p.\,143,
Cor.\,1}, or whenever $\dim R=2$, as follows e.g., from the Hoskin-Deligne
formula, see \cite{L3, p.\,223, (3.1.1)}.

\proclaim{Assumption (2.1)}
There exists a map $f\:Y\to\Spec(R)$ which factors as a sequence of
blowups with nonsingular centers, such that $I\O_Y$ is invertible.
\endproclaim

The basic conjecture ~(1.6) will be deduced from the following vanishing
conjecture.

\proclaim{Vanishing Conjecture (2.2)}
With $I$ and $f\:Y\to\Spec(R)$ as above,
$$
H^i(Y, I\omega_Y)=0\quad \text{for all }\, i>0.
$$
\endproclaim

\remark\nofrills{Remarks{\bf(2.2.1).}\usualspace} (a) Cutkosky has proved the
vanishing conjecture for
local rings essentially of finite type over a field of characteristic zero, see
\cite{C}. He uses Kodaira vanishing, which fails, in general, in positive
characteristic---but that does not preclude the conjecture holding for
special maps such as $f$.

(b) It was noted in~(1.3) that $\omega_Y$ is a dualizing sheaf for~$f$.
By duality \cite{L2, p.\,188}, the conjecture is equivalent to the vanishing of
$H^i_E(Y,(I\O_Y)^{-1})$ for all $i<\dim R$, where $E\set f^{-1}\{m\}$ is the
closed fiber. For $d=2$, this dual assertion is proved in
\cite{L2, p.\,177, Thm.\,2.4}.

(c) For $I=R$, the conjecture is a form of Grauert-Riemenschneider vanishing,
and is readily proved by induction on the number of blowups making up
the map~$f$. For arbitrary~$I$, the conjecture is equivalent to the
vanishing of~$H^i(Y,\Cal Q\mkern1mu)$ for all~$i>\nmb0$ and every invertible
quotient~$\Cal Q$ of a finite direct sum of copies of~$\omega_Y$
(because $I\O_Y$ is a quotient of a direct sum of copies of~$\O_Y$\dots)

Moreover, if $g\:Z\to\Spec(R)$ is the normalized blowup of~$I$ and
$h\:Y\to Z$ is the domination map, then using the Leray spectral sequence
for $f=gh$, and ampleness of~$I\O_Z$, one shows that the vanishing of
$R^ih_*\omega_Y\ (i>0)$
is equivalent to the vanishing of $H^i(Y, I^n\omega_Y)$ for all $n\gg0$.
In other words, Conjecture~(2.1) is somewhat stronger than
Grauert-Riemenschneider vanishing for ``sandwiched singularities."

(d) Theorem (4.1) of \cite{L6, p.\, 153} shows that there is
an $R$-ideal $L$ such that $Y$ in~(2.1) is the blowup of~$L$,
i.e., the Proj of the Rees ring $R[Lt]$ ($t$ an indeterminate), and such that
furthermore $R[Lt]$ is Cohen-Macaulay (CM). This leads to another conjecture
which can be shown to imply the vanishing one:

\proclaim{CM Conjecture}
Let $L=II'$, with $L$, $I$ and $I'$ integrally closed $R$-ideals,
and assume that $R[Lt]$ is CM and normal. Then for some~$e>0$, the
ideal $IR[L^et]$ $(\!$which is divisorial\/$\mkern1mu)$ is CM
as an $R[L^et]$-module.
\endproclaim

\nextpart{2.3} We show next that Conjecture (1.6) follows from the
vanishing conjecture.
\footnote{All we'll need here is that $Y$ is regular and $I\O_Y$ is
invertible.%
}
{\it Thus (1.6) does hold for local rings of
smooth points of algebraic  varieties in characteristic zero,
or when $\dim R=2$.} (See the preceding remarks (a) and~(b).)

We first reduce to the case where $R/m$ is infinite by passing, as usual,
to $S\set\nmb R[t]_{mR[t]}$. We have already seen, in proving (1.4.1)(iii),
that for any $R$-ideal ~$L$, $\tl LS\subset \tl{LS}$; but now we need equality,
which we get by applying Prop.\,(1.3.1) to $Y\otimes_R S$, with $Y$ as
in~(2.1).
(I don't know a more elementary way!)

Now let $I_0=(a_1,\dots,a_\ell)R$ be a reduction of~$I$, so that
$I_0\O_Y=I\O_Y$.
Let $F$ be the direct sum of $\ell$ copies of~$(I_0\O_Y)^{-1}$, and let
$\sigma\:F\to \O_Y$ be the $\O_Y$-homomorphism defined by the sequence~
$(a_1,\dots,a_\ell)$. Then we have a Koszul complex
$$
K(F,\sigma): \ 0\to\Lambda^\ell F\to\Lambda^{\ell-1} F \to \dots \to
\Lambda^1 F @>\sigma >> \O_Y \to 0
$$
(see \cite{LT, p.\,111}) which is locally split, so that
$K(F,\sigma)\otimes I^{n+1}\omega_Y\ (n\ge\ell-1)$ is exact.
By (2.2), and with $H^i(-)\set H^i(Y,-)$,
$$
H^1(I^{n-1}\omega_Y)=H^2(I^{n-2}\omega_Y)=\dots
 =H^{\ell-1}(I^{n+1-\ell}\omega_Y)=0.
$$
Hence, as in~\cite{LT, p.\,112, Lemma~(5.1)} we can conclude that
$$
H^0(I^{n+1}\omega_Y)=IH^0(I^n\omega_Y),
$$
i.e., by Proposition (1.3.1), $\tl{I^{n+1}}=I@,@,\tl{I^n}$.\qed
\smallskip

\subheading{3. Dimension 2}
Except in Lemma (3.2.1), $(R,m)$ will be a two-dimensional regular local ring
and $I$ will be an $m$-primary $R$-ideal.
The purpose of this section is to give a
number of alternative descriptions of~$\tl I$.

\nextpart{3.1}
It is pointed out in the footnote on p.\,235
of~\cite{L4} that when $I$ is a simple integrally closed ideal, the definition
of the adjoint of~$I$ given in \cite{L4, p.\,229} and \cite{L5, p.\, 299}
agrees with the one in this paper (see Proposition~(1.3.1)). Let us extend
this result---more specifically, the not-quite-correctly stated
Corollary~(4.1) of~\cite{L4, p.\,233}---to arbitrary $I$.

The {\it point basis\/} of~$I$ is the family of
integers~$(\ord_S(I^S))_{S\supset R}$ where $S$ runs through all
two-dimensional regular local rings between ~$R$ and its fraction field,
and $I^S\set\nmb\bigl(\gcd(IS)\bigr)^{-1}IS$, the {\it $S$-transform\/} of~$I$.
There are only finitely many~$S$ for
which $\ord_S(I^S)\ne 0$; these are called the {\it base points\/}
of~$I$ \cite {L4, p.\,225}. Two $m$-primary ideals $I'$ and $I''$ have the
same point basis iff their integral closures coincide \cite{L3, p.209, (1.10)}.

Consider a sequence of regular schemes
$$
\Spec(R)=:X_0 @<<f_0< X_1 @<<f_1< \cdots @<<f_n< X_{n+1}=:X
$$
where $f_i\:X_{i+1}\to X_i\ (0\le i\le n)$ is obtained by blowing up
a point on~$X_i$ whose local ring~$(S_i,m_i)$ is a base point of~$I$, and where
$I\O_X$ is invertible. Denote by~$m_i\O_X$ the invertible $\O_X$-ideal whose
stalk at $x\in X$ is $m_i\O_{X\!,x}$ if $\O_{X\!,x}\supset S_i\mkern1mu$,
and~$\O_{X\!,x}$ otherwise.  Then
$$
\omega_X^{-1} = \prod_{i=0}^n m_i\O_X\mkern1mu,
$$
see the end of the proof of Corollary~(1.4.4), or the footnote
in~\cite{L4, p.\,235}. Let $E_i\cong\nmb\Bbb P^1_{S_i/m_i}$
be the curve on~$X$ corresponding to the $m_i$-adic valuation;
and let $[S_i\:\mkern-4.5mu R]$ be the degree of
the field extension ~$(S_i/m_i)/(R/m)$.
The intersection number~$(E_i\cdot E_i)$ is $-d_i[S_i\:\mkern-4.5mu R]$
for some positive integer~$d_i\mkern1mu$,
and $d_i=1$ iff $I^{S_i}$~generates an invertible ideal on~$X_{i+1}$, i.e.,
iff $I^{S_i}$ is of the form $m_i^d\ (d>0)$, in which case
$(I\O_X\cdot E_i)=d\mkern1mu[S_i\:\mkern-4.5mu R]$. Moreover, as in \cite{L4,
p.\,235},
$$
(\omega_X\cdot E_i)=-(E_i\cdot E_i) -2[S_i\:\mkern-4.5mu R].
$$
It follows that $(I\omega_X\cdot E_i)\ge 0$ for all~$i$, and hence,
by \cite{L1, p.\,220, Thm.\,(12.1)(ii)},
{\it $I\omega_X$ is generated by its global sections,} i.e., by~$\tl I$, see
Prop.\,(1.3.1). Thus:
$$
I\omega_X= \tl I\O_X\,\quad\text{i.e.,}\quad
 I\O_X=\tl I\,\prod_{i=0}^n m_i\O_X\mkern1mu. \tag 3.1.1
$$

For any $S\supset R$, we have then
$$
\ord_S(I\mkern1mu)-\ord_S(\tl I\mkern1mu)=\sum_{S_i\subset S} \ord_S(m_i).
$$
On the other hand, setting, for any two-dimensional regular local~$T$
between $R$ and its fraction field, $r_T\set\ord_T(I^T)$,
$\mkern-.5mu\widetilde {\mkern.5mu r}_T\set\ord_T(\tl I\mkern1.5mu^T)$,  we
have
$$
\ord_S(I\mkern1mu)-\ord_S(\tl I\mkern1mu)=
\sum_{T\subset S} \ord_S(m_i)(r_T-\mkern-.5mu\widetilde {\mkern.5mu r}_T),
$$
see \cite{L4, p.\,301, Remark (1)}. By induction on the length of the unique
sequence of quadratic transforms
$R\set R_0\subset R_1\subset\dots\subset S$ (see \cite{A, p.\,343, Thm.\,3}),
we deduce that
$$
\spreadlines{1\jot}
\align
r_S-\mkern-.5mu\widetilde {\mkern.5mu r}_S &=1\qquad(S=S_1,S_2,\dots,S_n) \\
&=0\qquad\text{otherwise.}
\endalign
$$
But since $S_1,S_2,\dots,S_n$ are precisely the base points of~$I$, i.e.,
those $S$ such that $r_S>0$, what this amounts to is that
$\mkern-.5mu\widetilde {\mkern.5mu r}_S=\bigl(\max(0,\,r_S-1)\bigr)$.
Thus:

\proclaim{Proposition (3.1.2)}
$\tl I$ is the unique integrally closed ideal whose point basis is
$$
\bigl(\max(0,\,\ord_S(I^S)-1)\bigr)_{S\supset R}\,.
$$
\endproclaim

For any two-dimensional regular local~$T$
between $R$ and its fraction field, the point basis of the transform~$I^T$
is obtained from that of~$I$ by restriction to
those $S$ which contain~$T$. Moreover, a theorem of~Zariski states that
$I^T$ is integrally closed if $I$ is (see e.g., \cite{L5, p.\,300}). We
have then the following generalization of~\cite{L4, p.\,231, Thm.\,(3.1)}:

\proclaim{Corollary (3.1.3)}
Adjoint commutes with transform: for all~$T,$\ $\tl{I^T}=\nmb\tl I^T$.
\endproclaim

\nextpart{3.2} For the next result, let $I$ be any non-zero integrally closed
ideal in a $d$-dimensional
regular local ring~$R$, such that $I$ has a reduction~$I_0$
generated by a {\it regular sequence\/} $(a_1,a_2,\dots,a_d)$. Let $Y_0$ be the
blowup up of $I_0$, let $\pi\:Y\to Y_0$ be the normalization map, and let
$\C$ be the conductor of $Y$ in ~$Y_0$. Then $\C$ is independent of~$I_0$:

\proclaim{Lemma (3.2.1)}
With the preceding notation, we have\/
$\mkern1mu \C=I^{d-1}\omega_Y\,$.
\endproclaim

\proof
Noting that $Y_0\to\Spec(R)$ is a local complete intersection map (see footnote
under (1.3.2)(c)), and arguing as on pp.\;205--207 of~\cite{LS}, we find that
$$
\pi_*\omega_Y=\sHom(\pi_*\O_Y,\omega_{Y_0})=
\sHom(\pi_*\O_Y,(I_0\O_{Y_0})^{1-d}),
$$
so that
$$
\pi_*I^{d-1}\omega_Y=\pi_*I_0^{d-1}\omega_Y=I_0^{d-1}\pi_*\omega_Y
=\sHom(\pi_*\O_Y,\O_{Y_0})=\pi_*\C,
$$
whence the assertion.
\endproof

More can be said in the two-dimensional case.

\proclaim{Proposition (3.2.2)} With the preceding notation, when
$d=\dim R=2,$\ $\C$ is generated by its global sections
$H^0(Y,\C)=H^0(Y,I\omega_Y)=\tl I,$\ i.e.,
$\C=\tl I\O_Y$.
\endproclaim

\proof
Choose $X$ as in~(3.1), and let $g\:X\to Y$ be the domination map (which exists
because $I\O_X$ is invertible). As in the proof of~Prop.\,(1.3.1),
$g_*(I\omega_X)=I\omega_Y=\C$, the last equality by Lemma~(3.2.1).
Also, by \cite{L1, p.\,209, Prop.\,(6.5)}, the $\O_Y$-ideal $\tl I\O_Y$ is
integrally closed. Hence, and by ~(3.1.1),
$$
\tl I\O_Y=g_*(\tl I\O_X)=g_*(I\omega_X)=I\omega_Y=\C,
$$
whence the assertion.
\endproof
\medskip
Here is another characterization of $\tl I$.

\proclaim{Proposition (3.3)}
Let $I$ be an $m$-primary integrally closed ideal in a regular local ring~$R$
of dimension\/ $2,$\ let\/ $f\:Y\to\Spec(R)$
be the normalized blowup of\/~$I,$\ and let
$I_0=(a,b)R$ be a reduction of~$I$. Then with $D$ an injective hull
of\/~$R/m,$\ we have a duality isomorphism
$$
R/\tl I\mkern1mu\cong\Hom_R(I/I_0, D).
$$
Hence the $R$-module~$I/I_0$ depends only on~$I,$\ and its annihilator
$I_0\:\!I$ is just $\tl I$.
\endproclaim

\proof
Recall that $H^1(Y,\O_Y)=0$ \cite{L1, p.\,199, Prop.\,(1.2)}.
With  $\I\set I\O_Y=(a,b)\O_Y$, we have the exact
Koszul complex
$$
0\longrightarrow \I^{-1}@>-b\oplus a>>\O_Y\oplus\O_Y@>(a,\mkern1mu b)>>
\I\longrightarrow  0,
$$
whence an exact homology sequence, with $H^\bullet(-)\set H^\bullet(Y,-)$,
$$
R\oplus R=H^0(\O_Y\oplus\O_Y)@>(a,\mkern1mu b)>>H^0(\I)=I\to H^1(\I^{-1})\to
H^1(\O_Y\oplus\O_Y)=0,
$$
yielding $$I/I_0\cong H^1(\I^{-1}).$$
We already noted that $H^1(\O_Y)=0$, and since $f$ has fibers of dimension
$<2$ therefore $H^2(O_Y)=0$; so
$$H^1(\I^{-1})\cong H^1(\I^{-1}/\O_Y).$$
Further, with $E\set\nmb Y\otimes_R(R/m)$
the closed fiber,
we have that $\I^{-1}/\O_Y$ vanishes on~
$U\set\nmb Y\setminus\nmb E\cong\nmb\Spec(R)\setminus\{m\}$.
We conclude that
$$
H_E^1(\I^{-1}/\O_Y)\cong H^1(\I^{-1}/\O_Y)\cong I/I_0.
$$

Denoting the dualizing functor~$\Hom_R(-,D)$ by $-'$,  we have,
by \cite{L2, p.\,188},
$$
H^2_E(\I^{-1})\cong \text{Ext}^0(\I^{-1},\,\omega_Y)'\cong
H^0(\I\otimes\omega_Y)'
=(\tl I\,)\mkern1mu',
$$
and similarly
$$
H^2_E(\O_Y)\cong \text{Ext}^0(\O_Y,\,\omega_Y)'\cong H^0(\omega_Y)'=R'.
$$
Recall from (2.2.1)(b) that $H^1_E(\I^{-1})=0$. So there is an exact sequence
$$
0\to H_E^1(\I^{-1}/\O_Y)\to H^2_E(\O_Y)\to H^2_E(\I^{-1})\to 0
$$
whose dual is an exact sequence
$$
0\to \tl I\to R \to \Hom_R(I/I_0,D)\to 0
$$
which gives the desired conclusion.
\endproof

\end